\input amstex
\documentstyle{amsppt}
\TagsOnRight
\headline={\tenrm\hss\folio\hss}
\magnification=\magstep2
\vsize 23 true cm
\hsize 16 true cm

\centerline{\bf ON THE CONJUGACY SEPARABILITY } \centerline{\bf OF GENERALIZED FREE
PRODUCTS OF GROUPS}
\bigskip

\centerline{\bf E.~A.~Ivanova}
\bigskip

\smallskip

\hangindent 28 pt \hangafter -6 {\it It is proved that generalized free product of two finite
$p$-groups is a conjugacy $p$-separable group if and only if it is residually finite
$p$-groups. This result is then applied to establish  some sufficient conditions for conjugacy
$p$-separability of generalized free product of infinite groups.}
\bigskip

1. A group $G$ is called conjugacy separable (conjugacy $p$-separable) group if whenever
elements $a$ and $b$ of $G$ are not conjugate in $G$, there is a homomorphism $\varphi$ of $G$
onto finite (respectively, finite $p$-group) $X$ such that elements $a\varphi$ and $b\varphi$
are not conjugate in $X$.

It is easy to see that a conjugacy separable group is residually
finite group and a conjugacy $p$-separable group is residually
finite $p$-groups. In general inverse statements are not true. But
in some cases residually finite group (residually finite
$p$-groups) is a conjugacy separable (respectively, conjugacy
$p$-separable) group too. For example, J.~Dyer [7] have proved
that a free product with amalgamated subgroups of two finite
groups is conjugacy separable group (in 1963 G.~Baumslag [5]
showed that any such group is residually finite).

Not every free product with amalgamated subgroups of two finite $p$-groups must be residually
finite $p$-groups. G.Higman [9] have obtained necessary and sufficient conditions for such
groups to be residually finite $p$-groups. The question arises whether these conditions are
enough for such group to be a conjugacy $p$-separable group.

An element $g$ of a group $G$ is called $C_{fp}$-separable if for
any $a\in G$ such that elements $a$ and $g$ are not conjugate in
$G$, there exists a homomorphism $\varphi$ of $G$ onto finite
$p$-group $X$ such that $a\varphi$ and $g\varphi$ are not
conjugate in $X$. So a group $G$ is conjugacy $p$-separable if and
only if every element $g\in G$ is $C_{fp}$-separable.

It was proved in [1] that
\newline {\it if a free product with amalgamated subgroups of two finite
$p$-groups is a residually finite $p$-group then every infinite order element $g\in G$ is
$C_{fp}$@-separable.}

In fact the following generalization of this statement holds:

\proclaim{\indent Theorem 1} Suppose $G=\bigl(H
* K; \ A=B,\ \varphi\bigr)$ is a free product of two finite $p$-groups $H$ and $K$ with
amalgamated via isomorphism $\varphi$ subgroups $A$ and $B$. $G$ is a conjugacy $p$-separable
group if and only if $G$ is a residually finite $p$-groups.
\endproclaim

Applying this result and using a standard technique we proved the
following theorem:

\proclaim{\indent Theorem 2} Suppose that $H$ and $K$ are conjugacy $p$-separable groups,
$A\leqslant H$ and $B\leqslant K$ are central subgroups and every finite $p$-index subgroup of
$A$ and of $B$ is $p$-separable in $H$ and in $K$ respectively. Then $G=\bigl(H * K; \ A=B, \
\varphi\bigr)$ is a conjugacy $p$-separable group.
\endproclaim

(Recall that a subset $M$ of a group $G$ is called $p$-separable if for every $a\in G$,
$a\notin M$, there is a homomorphism $\varphi$ of $G$ onto finite $p$-group $X$ such that
$a\varphi\notin M\varphi$.)\smallskip

The description of conjugacy $p$-separable finitely generated nilpotent groups has been given
in [2]. It follows from the theorem 2 that the following statement holds:

\proclaim{\indent Theorem 3} Suppose $G=\bigl(H
* K; \ A=B,\ \varphi\bigr)$ is a free product with amalgamated subgroups
of two finitely generated nilpotent groups $H$ and $K$, $H$ and $K$ are conjugacy
$p$-separable groups, $A$ and $B$ are $p^{\prime}$-isolated central subgroups of $H$ and $K$
respectively. Then $G$ is a conjugacy $p$@-separable group.
\endproclaim

(Recall that if $p$ is a prime then a subgroup $X$ of a group $Y$ is called $p$-isolated if
$y^p\in X$ implies $y\in X$ for every $y\in Y$. A subgroup $X$ of a group $Y$ is called
$p^{\prime}$-isolated if $X$ is a $q$-isolated for every prime $q\neq p$.)
\bigskip

2. The  proof of theorem 1 is a certain modification of the J.~Dyer's proof of her result in
[8].

Graph $\Gamma$ is a system of two sets $V=V(\Gamma)$ (the set of vertexes) and $E=E(\Gamma)$
(the set of edges) and of three mappings $\overline{ }:E\to E$, $o:E\to V$ and $t:E\to V$ such
that $o(\overline{e})=t(e)$, $t(\overline{e})=o(e)$, $\overline{e}\neq e$ and
$\overline{\overline{e}}=e$ for every $e\in E$. The edge $\overline{e}$ is called inverse to
edge $e$, the vertex $o(e)\in V$ is called the origin of edge $e$, the vertex $t(e)\in V$ is
called the end of edge $e\in E$. If $o(e)=u$ and $t(e)=v$ then we write $e=(u,v)$.

The group graph is a pair $(\Cal G, \Gamma)$ of connected graph
$\Gamma$ and mapping $\Cal G$. The mapping $\Cal G$ associates
every vertex $v\in V(\Gamma)$ with group $G_v$ and every edge
$e\in E(\Gamma)$, $e=(u,v)$ with group $G_e$ and two mappings
$\rho_e : G_e\to G_u$ and $\tau_e : G_e\to G_v$ such that
$G_{\overline{e}}=G_e$, $\rho_{\overline e} =\tau_e$ and
$\tau_{\overline e} =\rho_e$. The groups $G_v$ and $G_e$ are
called vertex group and edge group of a group graph $(\Cal G,
\Gamma)$ respectively.

Suppose $(\Cal G, \Gamma)$ is a group graph and $T$ is a maximal
tree of $\Gamma$. Let $X_v$ be a set of generators and $R_v$ be a
set of relations of the vertex group $G_v$ for every vertex $v\in
V=V(\Gamma)$ (and if $v_1\neq v_2$ then $X_{v_1}\cap X_{v_2}=
\emptyset$). The fundamental group $\pi(\Cal G, \Gamma)$ of a
group graph $(\Cal G, \Gamma)$ is a group with generators
$\bigcup_{v\in V}X_v$ and $t_e$, $e\in E(\Gamma)\setminus E(T)$,
and relations $\bigcup_{v\in V}R_v$ and
$$
\aligned g=g(\rho_e^{-1}\tau_e) &, \qquad e\in E(T),\ g\in G_e\rho_e,\\
t_e^{-1}gt_e=g(\rho_e^{-1}\tau_e) &, \qquad e\in E(\Gamma)\setminus E(T),\ g\in G_e\rho_e, \\
t_{_{\overline{e}}}=t_e^{-1} &, \qquad e\in E(\Gamma)\setminus
E(T).
\endaligned
$$
It is possible to prove that the group $\pi(\Cal G, \Gamma)$ does
not depend on a choice of vertex groups presentations and maximal
tree $T$.

It is well known ([6, 10, 11]) that every finite extension of a
free group is isomorphic to a fundamental group of a group graph
with finite vertex groups.

To prove the theorem 1 the following result is also required
([1]):

\proclaim{\indent Proposition 2.1} Suppose $H\leqslant G$ is a subnormal finite $p$-index
subgroup. If $h\in H$ is a $C_{fp}$-separable in group $H$ then $h$ is $C_{fp}$-separable in
group $G$.
\endproclaim

Now we are ready to prove the theorem 1.

Suppose $G=\bigl(H * K; \ A=B,\ \varphi\bigr)$ is a free product
of two finite $p$@-groups $H$ and $K$ with amalgamated via
isomorphism $\varphi$ subgroups $A$ and $B$. The necessary
conditions in theorem are evident.

Let $G$ be a residually finite $p$-groups. Then ([4, lemma 2.1])
$G$ is an extension of a free group $F$ by finite $p$-group. In
view of mentioned above result from [1] to prove that $G$ is a
conjugacy $p$-separable it is enough to show that whenever $a$ and
$b$ are finite order elements and not conjugate in $G$, there is a
homomorphism $\psi$ of $G$ onto a finite $p$-group $X$ such that
$a\varphi$ and $b\varphi$ are not conjugate in $X$.

Since subgroup of $G$ which is generated by $F$ and $a$ is a subnormal in $G$ it follows from
proposition 2.1 that we can assume that group $G$ is generated by $F$ and $a$. Then the
quotient group $G/F$ is cyclic and therefore if $aF\neq bF$ then natural homomorphism of $G$
onto $G/F$ is required.

Suppose $aF=bF$. Since $F$ is a torsion-free so the orders of $a$
and $b$ are equal to $p^n$ ($n\geqslant 1$). By remark above the
group $G$ is isomorphic to a fundamental group of a group graph
and its vertex subgroups can be embedded into the cyclic group of
order $p^n$. D.~Dayer [8] showed that in this case there is a
homomorphism $\psi$ of $G$ onto fundamental group $H=\pi(\Cal H,
\Gamma)$ of a group graph $(\Cal H, \Gamma)$ such that

1) the graph $\Gamma$ has only two vertexes $u$ and $v$;

2) the vertex groups $H_u$ and $H_v$ are cyclic of order $p^n$ and they are
generated by elements $x=a\varphi$ and $y=b\varphi$ respectively;

3) the order of every edge group $H_e$ is less than $p^n$.

Thus the generators of $H$ are $x$, $y$ and $t_e$, where edge
$e\in E(\Gamma)$ is not equal to some fixed edge, the relations of
$H$ are:

a) $x^{p^n}=1$, $y^{p^n}=1$;

b) $x^r=y^s$, where elements $x^r\in H_u$ and $y^s\in H_v$ have the same order which is less
than $p^n$;

c) $t_e^{-1}h_1t_e=h_2$, where the elements $h_1\in H_u$ and $h_2\in H_v$ have the same order
which is less than $p^n$.

By condition 3) every element $x^r, h_1\in H_u$, $y^s, h_2\in H_v$
from relations b) and c) belongs to subgroup $K_u$ or $K_v$
respectively, where orders of $K_u$ and $K_v$ are equal to $p^k$,
$k<n$. Therefore relations b) and c) are trivial in $H$ modulo $N$
where $N$ is a normal closure of $K_u$ and $K_v$. Thus the
quotient-group $H/N$ is a free product of finite cyclic groups
$H_u/K_u=(x)$ and $H_v/K_v=(y)$, which orders are equal to
$p^{n-k}$, and the set of infinite cyclic groups $(t_e)$. The
composition of $\psi$, natural homomorphism of $H$ onto $H/N$ and
an evident homomorphism of $H/N$ onto direct product of the groups
$H_u/K_u$ and $H_v/K_v$ is required homomorphism of $G$. The
theorem 1 is proved.
\bigskip

3. Suppose $H$ and $K$ are groups, $A$ is a subgroup of $H$, $B$
is a subgroup of $K$, $\varphi:A\to B$ is a isomorphism.

Every element $x$ from $G=(H*K;\ A=B,\ \varphi)$ can be presented as $x=x_1x_2\cdot\dots\cdot
x_n$, where each $x_1$, $x_2$,\dots , $x_n$ is from factor $H$ or $K$ and if $n>1$ then  $x_i$
and $x_{i+1}$ are from different factors for every $i=1,\dots,n-1$ (therefore they do not
belong to $A$ and $B$). This presentation is called reduced form of $x$ and the number $n$
(that doesn't depend of choice of such presentation) is called the length of $x$. An element
$x\in G$ is called cyclically reduced if either its length $n$ equals to 1 or $n>1$ and
elements $x_1$ and $x_n$ in its reduced form are from different factors $H$ and $K$. In this
case the expression $u_i=x_ix_{i+1}\cdot\dots\cdot x_nx_1\cdots x_{i-1}$ is reduced for each
$i=1, 2,\dots, n$. The element $u_i$ is called a cyclic permutation of $x$ (if $n=1$ then $x$
is a unique cyclic permutation of $x$).

When amalgamating subgroups $A$ and $B$ are central then the general conditions for two
elements from $G$ to be conjugate ([3]) can be simplified:

\proclaim{\indent Proposition 3.1} Suppose $G=(H*K;\ A=B, \varphi)$ is a free product of two
groups with amalgamated central subgroups $A$ and $B$. For every $g\in G$ there is a
cyclically reduced $x\in G$ such that $g$ and $x$ are conjugate. Suppose $x\in G$ and $y\in G$
are cyclically reduced. Then $x$ and $y$ are conjugate in $G$ if and only if their length are
equal and either they are from one factor $H$ or $K$ and conjugate in it, or their lengths
more then 1 and one of these elements equals to the cyclic permutation of another.
\endproclaim

Let's remind also the following notion ([5]). Subgroups $R\leqslant H$ and $S\leqslant
K$ are called $(A,B,\varphi)$-compatible if $(A\cap R)\varphi =B\cap S$. If normal subgroups
$R\leqslant H$ and $S\leqslant K$ are $(A,B,\varphi)$-compatible then the mapping
$\varphi_{{}_{R,S}}: AR/R \to BS/S,$ where $(aR)\varphi_{{}_{R,S}}=(a\varphi )S$ ($a\in A$),
is an isomorphism of subgroup $AR/R\leqslant H/R$ onto subgroup $BS/S\leqslant K/S$. Thus
there is a free product
$$
G_{R,S} = \left(H/R * K/S; \  AR/R = BS/S, \ \varphi_{{}_{R,S}}\right)
$$
of the groups $H/R$ and $K/S$ with amalgamated via isomorphism
$\varphi_{{}_{R,S}}$ subgroups $AR/R$ and $BS/S$. Natural
homomorphisms of the group $H$ onto quotient group $H/R$ and of
the group $K$ onto quotient group $K/S$ can be extended to a
homomorphism $\rho_{{}_{R,S}}$ of the group $G=(H*K;\ A=B,\
\varphi)$ onto the group $G_{R,S}$.

\proclaim{\indent Proposition 3.2} Suppose that $H$ and $K$ are conjugacy $p$-separable
groups, $A\leqslant H$ and $B\leqslant K$ are central subgroups and subgroups $A$ and $B$ and
every finite $p$-index subgroup of $A$ and $B$ are $p$-separable in group $H$ and $K$
respectively. Then for every finite $p$-index normal subgroups $M\leqslant H$ and $N\leqslant
K$ there are $(A, B, \varphi )$-compatible finite $p$-index normal subgroups $R\leqslant H$
and $S\leqslant K$ such that $R\leqslant M$ and $S\leqslant N$.
\endproclaim

\demo{\indent Proof} Suppose $M$ contains a finite $p$-index (in $A$) subgroup $U\leqslant A$.
Then subgroup $U$ is $p$-separable in the group $H$ and therefore the quotient group $H/U$ is
residually finite $p$-groups. Also since the quotient group $A/U$ is a finite subgroup of the
quotient group $H/U$, therefore there is a finite $p$-index normal subgroup $R/U$ of the group
$H/U$ such that $R/U\cap A/U=1$. Then $R$ is a finite $p$-index normal subgroup of the group
$H$ and $R\cap A=U$. We can consider also that $R\leqslant M$. The similar reasoning is fair
and for the group $K$.

Suppose $M$ and $N$ are finite $p$-index normal subgroups of the
groups $H$ and $K$ respectively, $U=(M\cap A)\cap (N\cap
B)\varphi^{-1}$ and $V=(M\cap A)\varphi\cap (N\cap B)$. Then there
are finite $p$-index normal subgroups $H$ and $K$ of the groups
$R$ and $S$ respectively such that $R\leqslant M$, $R\cap A=U$,
$K$$S\leqslant N$ and $S\cap B=V$. Since $U\varphi = V$ the
subgroups $R$ and $S$ are required.
\enddemo

Now we are ready to prove the theorem 2. Suppose $H$ and $K$ are
conjugacy $p$-separable groups, $G = \bigl(H * K; \ A = B,
\,\varphi \bigr)$ is a free product of the groups $H$ and $K$ with
amalgamated central subgroups $A$ and $B$.

Using standard methods of the proof for the free product of two
groups with amalgamated subgroups to be residually finite
$p$-groups it is easy to receive the following statement:

\proclaim{\indent Proposition 3.3} Suppose $H$ and $K$ are
residually finite $p$-groups, $A\leqslant H$ and $B\leqslant K$
are central subgroups and subgroups $A$ and $B$ and every finite
$p$-index subgroup of $A$ and $B$ are $p$-separable in group $H$
and $K$ respectively. Then $G=\bigl(H\ *\ K;\ A=B,\ \varphi\bigr)$
is a residually finite $p$-groups.
\endproclaim

It follows from the theorem 1 that if $R\leqslant H$ and
$S\leqslant K$ are $(A,B,\varphi)$-compatible finite $p$-index
normal subgroups then $G_{R,S}$ is a conjugacy $p$-separable group
(it is proved in [9] that $G_{R,S}$ is a residually finite
$p$-groups). Therefore it is enough to prove that whenever $x\in
G$ and $y\in G$ are not conjugate in $G$ there are
$(A,B,\varphi)$-compatible finite $p$-index normal subgroups
$R\leqslant H$ and $S\leqslant K$ such that $x\rho_{{}_{R,S}}$ and
$y\rho_{{}_{R,S}}$ are not conjugate in $G_{R,S}$.

Suppose $f\in G$ and $g\in G$ are not conjugate in $G$. Since
proposition 3.1 we can assume without generality loss that $f$ and
$g$ are cyclically reduced. Let's consider some cases.

{\it Case 1.} The lengths of $f$ and $g$ are not equal.

Since the subgroups $A$ and $B$ are $p$-separable in the groups $H$ and $K$ respectively there
are finite $p$-index normal subgroups $M\leqslant H$ and $N\leqslant K$ such that all factors
in the reduced forms of $f$ and $g$ (it is fair only for the elements of length 1) are not
from $AM$ and $BN$ respectively. Since proposition 3.2 there are $(A, B, \varphi )$-compatible
finite $p$-index normal subgroups $R\leqslant H$ and $S\leqslant K$ such that $R\leqslant M$
and $S\leqslant N$. Then $f\rho_{{}_{R,S}}$ and $g\rho_{{}_{R,S}}$ are cyclically reduced in
the group $G_{R,S}$, the lengths of $f\rho_{{}_{R,S}}$ and $g\rho_{{}_{R,S}}$ are equal to the
lengths of $f$ and $g$ respectively and different. Therefore since proposition 3.1
$f\rho_{{}_{R,S}}$ and $g\rho_{{}_{R,S}}$ are not conjugate in $G_{R,S}$.

{\it Case 2.} The lengths of $f$ and $g$ are equal to 1, $f$ and
$g$ are from different factors $H$ and $K$.

Suppose $f\in H\setminus A$ and $g\in K\setminus B$. Since $A$ and
$B$ are $p$-separable subgroups in $H$ and $K$ respectively there
are finite $p$-index normal subgroups $M\leqslant H$ and
$N\leqslant K$ such that $f\notin AM$ and $g\notin BN$. At that
time since proposition 3.2 there are $(A,B,\varphi)$-compatible
finite $p$-index normal subgroups $R\leqslant H$ and $S\leqslant
K$ such that $R\leqslant M$ and $S\leqslant N$. Then
$f\rho_{{}_{R,S}}$ and $g\rho_{{}_{R,S}}$ are from different
factors $H/R$ and $K/S$ of $G_{R,S}$ and since proposition 3.1 the
elements $f\rho_{{}_{R,S}}$ and $g\rho_{{}_{R,S}}$ are not
conjugate in $G_{R,S}$.

{\it Case 3.} The lengths of $f$ and $g$ are equal to 1, $f$ and
$g$ are both from one factor $H$ or $K$.

Suppose $f\in H$ and $g\in H$. Since $f$ and $g$ are not conjugate
in $H$ and $H$ is a conjugacy $p$-separable group there is a
finite $p$-index normal subgroup $M\leqslant H$ such that $fM$ and
$gM$ are not conjugate in $H/M$.  At that time since proposition
3.2 there are $(A,B,\varphi)$-compatible finite $p$-index normal
subgroups $R\leqslant H$ and $S\leqslant K$ such that $R\leqslant
M$. Then $f\rho_{{}_{R,S}}$ and $g\rho_{{}_{R,S}}$ are both from
the factor $H/R$ of the group $G_{R,S}$ and the elements
$f\rho_{{}_{R,S}}$ and $g\rho_{{}_{R,S}}$ are not conjugate in
$H/R$. Therefore since proposition 3.1 $f\rho_{{}_{R,S}}$ and
$g\rho_{{}_{R,S}}$ are not conjugate in $G_{R,S}$.

{\it Case 4.} The lengths of $f$ and $g$ are equal and more then
1.

Suppose $g_1$, $g_2$, \dots , $g_r$ ($r$ is a length of $g$) are
all cyclic permutations of $g$. Since the subgroups $A$ and $B$
are $p$-separable in the groups $H$ and $K$ respectively there are
finite $p$-index normal subgroups $R_0\leqslant H$ and
$S_0\leqslant K$ such that all factors in the reduced forms of $f$
and $g$ are not from $AR_0$ and $BS_0$ respectively. Since $f$ is
not equal to $g_1$, $g_2$, \dots , $g_r$ and $G$ is a residually
finite $p$-groups (proposition 3.3) there is a finite $p$-index
normal subgroup $N\leqslant G$ such that $fN$ is not equal to
$g_1N$, $g_2N$, \dots , $g_rN$. Suppose $R=R_0\cap N$ and
$S=S_0\cap N$. Then $f\rho_{{}_{R,S}}$ and $g\rho_{{}_{R,S}}$ are
cyclically reduced in $G_{R,S}$, their lengths are equal to $r$
and $f\rho_{{}_{R,S}}$ is not equal to $g_1\rho_{{}_{R,S}}$,
$g_2\rho_{{}_{R,S}}$, \dots , $g_r\rho_{{}_{R,S}}$. Since every
cyclic permutation of $g\rho_{{}_{R,S}}$ is equal to one of the
elements $g_1\rho_{{}_{R,S}}$, $g_2\rho_{{}_{R,S}}$, \dots ,
$g_r\rho_{{}_{R,S}}$ the elements $f\rho_{{}_{R,S}}$ and
$g\rho_{{}_{R,S}}$ are not conjugate in $G_{R,S}$ (proposition
3.1).

The theorem 2 is proved.

Since every $p^{\prime}$-isolated subgroup of a finitely generated
nilpotent group is $p$-separable the theorem 3 immediately follows
from the theorem 2.

\bigskip
\centerline{\bf References}
\medskip

\noindent \item {1.} {\it Ivanova~E.~A.} On conjugacy
p-separability of free products  with amalgamated subgroups of two
groups, Math. notes, V.~76, (4), 2004, P.~502-509. (Russian)

\noindent \item {2.} {\it Ivanova~E.~A., Moldavanskii~D.~I.} On
conjugacy separability of finitely generated nilpotent groups,
Bull. Ivanovo State University, (3), 2004, P.~125-130. (Russian)

\noindent \item {3.} {\it Magnus~W.' Karrass~A., Solitar~D.}
Combinatorial group theory: Presentations of groups in terms of
generators and relations, Interscience publishers, New York,
London, Sydney, 1966.

\noindent \item {4.} {\it Moldavanskii~D.~I.} Finite residuality
of HNN-extensions, Bull. Ivanovo State University, (3), 2000'
P.~129-140. (Russian)

\noindent \item {5.} {\it Baumslag G.} On the residual finiteness
of generalized free pro\-ducts of nilpotent groups, Trans. Amer.
Math. Soc., V.106, 1963, P.~193-209.

\noindent \item {6.} {\it Cohen D.~E.} Groups with free subgroups
of finite index, Conf. Group Theory, Univ. Wiskonsin, Parkside,
1972, Lecture Notes Math., 1973, {\bf 319}, P.~26-44.

\noindent \item {7.} {\it Dyer J.~L.} Separating conjugates in
amalgamating free products and HNN-extensions, J\. Aust\. Math\.
Soc\., V. 29, (1), 1980, P.~35-51.

\noindent \item {8.} {\it Dyer J.~L.} Separating conjugates in
free-by-finite groups, J. London Math. Soc., V.20, (2), 1979,
P.~215-221.

\noindent \item {9.} {\it Higman G.} Amalgams of $p$@-groups, J.
Algebra, V.1, 1964, P.~301-305.

\noindent \item {10.} {\it Karrass A., Pietrowski A., Solitar D.}
Finite and infinite cyclic ex\-ten\-sions of free groups, J\.
Aust\. Math\. Soc., V. 16, 1973, P.~458-466.

\noindent \item {11.} {\it Scott~G.~P.} An embedding theorem for
groups with a free subgroup of finite index, Bull. Lond. Math.
Soc., V. 6, 1974, P.~304-306.
\medskip

Department of Mathematics, Ivanovo State University.

ea-ivanova\@mail.ru

\end